\documentclass[12pt, reqno]{amsart}
\usepackage{amsmath, amstext, amsbsy, amssymb, amscd}

\newtheorem{theorem}{Theorem}[section]
\newtheorem{lemma}[theorem]{Lemma}
\newtheorem{proposition}[theorem]{Proposition}
\newtheorem{corollary}[theorem]{Corollary}
\theoremstyle{definition}

\theoremstyle{remark}
\newtheorem{remark}[theorem]{Remark}

\newtheorem{notation}[theorem]{Notation}
\newtheorem{conjecture}[theorem]{Conjecture}

\numberwithin{equation}{section}

\newcommand{\C}{\mathbb C}
\newcommand{\ch}{\text{\rm ch}}

\newcommand{\id}{{\rm Id}}

\newcommand{\och}{\mathfrak{ch}}

\newcommand{\Q}{\mathbb Q}

\newcommand{\Supp}{{\rm Supp}}

\newcommand{\w}{\tilde}
\newcommand{\wtau}{\tilde \tau}
\newcommand{\W}{\widetilde}

\newcommand{\Z}{\mathbb Z}

{\vskip-\lastskip\medskip
  \noindent
  {\em #1.}\enspace
  }%
{\qed\par\medskip
  }

\begin{document}
\title[Gromov-Witten and Donaldson-Thomas correspondence] 
{The Gromov-Witten and Donaldson-Thomas correspondence
for trivial elliptic fibrations}

\author[Dan Edidin]{Dan Edidin$^1$}
\address{Department of Mathematics, University of Missouri, 
Columbia, MO 65211, USA} \email{edidin@math.missouri.edu}
\thanks{${}^1$Partially supported by an NSA grant}

\author[Zhenbo Qin]{Zhenbo Qin$^2$}
\address{Department of Mathematics, University of Missouri, 
Columbia, MO 65211, USA} \email{zq@math.missouri.edu}
\thanks{${}^2$Partially supported by an NSF grant}


\begin{abstract}
We study the Gromov-Witten and Donaldson-Thomas correspondence
conjectured in \cite{MNOP1, MNOP2} for trivial elliptic fibrations.
In particular, we verify the Gromov-Witten and Donaldson-Thomas 
correspondence for primary fields when the threefold is $E \times S$ 
where $E$ is a smooth elliptic curve and $S$ is a smooth surface 
with numerically trivial canonical class.
\end{abstract}

\maketitle
\date{}

\section{\bf Introduction}

The correspondence between the Gromov-Witten theory and
Donaldson-Thomas theory for threefolds was 
conjectured and studied in \cite{MNOP1, MNOP2}. 
Since then, it has been investigated extensively
(see \cite{MP, JLi, Kat, KLQ, Beh, BF2} and 
the references there). 
A relationship between the quantum cohomology of 
the Hilbert scheme of points in the complex plane and 
the Gromov-Witten and Donaldson-Thomas correspondence 
for local curves was proved in \cite{OP2, OP3}. 
The equivariant version was proposed and partially 
verified in \cite{BP, GS}. 
In this paper, we study the Gromov-Witten and 
Donaldson-Thomas correspondence when the threefold 
admits a trivial elliptic fibration.

To state our results, we introduce some notation
and refer to Subsect.~\ref{subsect_gw_inv} and
Subsect.~\ref{subsect_dt_inv} for details.
Let $X$ be a complex threefold, 
$\gamma_1, \ldots, \gamma_r \in H^*(X; \Q)$,
\begin{eqnarray*}  
\beta \in H_2(X; \Z) \backslash \{0\},
\end{eqnarray*}
$k_1, \ldots, k_r$ be nonnegative integers, 
and $u, q$ be formal variables. Let
\begin{eqnarray*}  
{\bf Z}_{\text{GW}}' \left (X; u|
\prod_{i=1}^r \tau_{k_i}(\gamma_i) \right )_\beta,
\qquad  
{\bf Z}_{\text{DT}}'\left (X; q| 
\prod_{i=1}^r \wtau_{k_i}(\gamma_i) \right )_\beta
\end{eqnarray*}
be the reduced degree-$\beta$ partition functions 
for the descendent Gromov-Witten invariants
and Donaldson-Thomas invariants of $X$ respectively.

\begin{conjecture}  \label{conj_pri}
(\cite{MNOP1, MNOP2}) Let $\beta \in H_2(X; \Z) 
\backslash \{0\}$ and $\mathfrak d 
= -\displaystyle{\int_{\beta} K_X}$.
Then after the change of variables $e^{iu} = -q$, we have
\begin{eqnarray}  \label{conj.1}
(-iu)^{\mathfrak d} \; {\bf Z}_{\text{GW}}'
\left (X; u|\prod_{i=1}^r \tau_{0}(\gamma_i) \right )_\beta 
= (-q)^{-\mathfrak d/2} \; {\bf Z}_{\text{DT}}'\left (X; q| 
\prod_{i=1}^r \wtau_{0}(\gamma_i) \right )_\beta.
\end{eqnarray}
\end{conjecture}

\begin{theorem} \label{thm}
Let $f: X = E \times S \to S$ be the projection 
where $E$ is an elliptic curve and $S$ is a smooth surface. 
Then the Gromov-Witten/Donaldson-Thomas 
correspondence (\ref{conj.1}) holds if either 
$\displaystyle{\int_{\beta} K_X = \int_{\beta} f^*K_S =0}$, 
or 
\begin{eqnarray*}  
\gamma_1, \ldots, \gamma_r \in f^*H^*(S; \Q) 
\subset H^*(X; \Q).
\end{eqnarray*}
\end{theorem}
\noindent
{\it Proof.}
The conclusion follows from Proposition~\ref{prop_gw2}
and Proposition~\ref{prop_dt2} when
\begin{eqnarray*}  
\int_{\beta} K_X = \int_{\beta} f^*K_S =0.
\end{eqnarray*}
It follows from Proposition~\ref{prop_gw1}
and Proposition~\ref{prop_dt1} when
\begin{equation}
\gamma_1, \ldots, \gamma_r \in f^*H^*(S; \Q) 
\subset H^*(X; \Q).                  
\tag*{$\qed$}
\end{equation}

\begin{corollary} \label{cor}
Let $E$ be an elliptic curve and $S$ be a smooth surface
with numerically trivial canonical class $K_S$. 
Then the Gromov-Witten/Donaldson-Thomas 
correspondence (\ref{conj.1}) holds for the threefold
$X = E \times S$.
\end{corollary}

In fact, when $\gamma_1, \ldots, \gamma_r \in f^*H^*(S; \Q) 
\subset H^*(X; \Q)$, Proposition~\ref{prop_gw1}
and Proposition~\ref{prop_dt1} state that after 
the change of variables $e^{iu} = -q$, 
\begin{eqnarray*}  \label{conj_des}
(-iu)^{\mathfrak d - \sum_i k_i} \; {\bf Z}_{\text{GW}}'
\left (X; u|\prod_{i=1}^r \tau_{k_i}(\gamma_i) \right )_\beta 
= (-q)^{-\mathfrak d/2} \; {\bf Z}_{\text{DT}}'\left (X; q| 
\prod_{i=1}^r \wtau_{k_i}(\gamma_i) \right )_\beta.
\end{eqnarray*}
This is consistent with (and partially sharpens) 
the Conjecture~4 in \cite{MNOP2} which is about 
the Gromov-Witten and Donaldson-Thomas correspondence 
for descendent fields. It would be interesting to see
whether this sharpened version holds for general 
cohomology classes $\gamma_1, \ldots, \gamma_r \in 
H^*(X; \Q)$.

Proposition~\ref{prop_gw1} and Proposition~\ref{prop_dt1}
are proved in Sect.~\ref{sect_gw} and Sect.~\ref{sect_dt}
respectively. The idea is to view the elliptic curve $E$ 
as an algebraic group and to use the action of $E$ on 
the moduli space $\overline{\mathfrak M}_{g, r}(X, \beta)$ 
of stable maps and the moduli space 
$\mathfrak I_n(X, \beta)$ of ideal sheaves.
The $E$-action on $\overline{\mathfrak M}_{g, r}(X, \beta)$
has no fixed points when $r \ge 1$, or $g \ne 1$, 
or $\beta \ne d \beta_0$. It follows from 
Lemma~\ref{lma_gw=0} that the corresponding Gromov-Witten
invariants are zero. The only exception is 
$\langle \rangle_{1, d \beta_0}$ which can be computed 
directly by using the work of Okounkov-Pandharipande \cite{OP1}
on the Gromov-Witten invariants of an elliptic curve 
and G\" ottsche's formula for the Euler characteristics of
the Hilbert scheme $S^{[d]}$ of points on a smooth surface $S$.
Similarly, the $E$-action on $\mathfrak I_n(X, \beta)$
has no fixed points when $n \ge 1$ or $\beta \ne d \beta_0$. 
It follows from Lemma~\ref{lma_dt_0} that  
the corresponding Donaldson-Thomas invariants are 
also zero. The only exception is 
$\langle \rangle_{0, d \beta_0}$ which can be computed 
directly by determining the obstruction bundle over 
the moduli space $\mathfrak I_0(X, d \beta_0) \cong S^{[d]}$.

It is expected that our approach can be used to handle
the {\it relative} Gromov-Witten and Donaldson-Thomas 
correspondence (see \cite{MNOP2}) for trivial elliptic 
fibrations. In another direction, one might attempt to 
study the (absolute and relative) Gromov-Witten and 
Donaldson-Thomas correspondence for nontrivial elliptic 
fibrations. We leave these to the interested readers.

\bigskip\noindent
{\bf Acknowledgments.} The authors thank Professors 
Robert Friedman, Yuan-Pin Lee, and Wei-Ping Li for 
valuable help.

\section{\bf Gromov-Witten theory} 
\label{sect_gw}
\subsection{\bf Gromov-Witten invariants} 
\label{subsect_gw_inv}
\par
$\,$

Let $X$ be a smooth projective complex variety. 
Fix $\beta \in H_2(X; \Z)$. 
Let $\overline{\mathfrak M}_{g, r}(X, \beta)$
be the moduli space of stable maps from connected genus-$g$
curves with $r$ marked points to $X$ representing 
the class $\beta$. The virtual fundamental class 
$[\overline{\mathfrak M}_{g, r}(X, \beta)]^{\rm vir}$
has been constructed in \cite{BF1, LT}. By ignoring the extra 
notation of stacks, the virtual fundamental class 
$[\overline{\mathfrak M}_{g, r}(X, \beta)]^{\rm vir}$ 
is defined by the element
\begin{eqnarray}  \label{pot_gw}
R(\pi_{g,r})_*({\rm ev}_{r+1})^*T_X
\end{eqnarray}
in the derived category $\mathfrak D_{\rm coh}
(\overline{\mathfrak M}_{g, r}(X, \beta))$ of coherent 
sheaves on $\overline{\mathfrak M}_{g, r}(X, \beta)$, where 
\begin{eqnarray*}
{\rm ev}_i: \; \overline{\mathfrak M}_{g, r+1}(X, \beta) 
\to X
\end{eqnarray*}
is the $i$-th evaluation map, and $\pi_{g,r}$ stands for 
the morphism:
\begin{eqnarray}   \label{pigr}
\pi_{g,r}: \; \overline{\mathfrak M}_{g, r+1}(X, \beta) 
\to \overline{\mathfrak M}_{g, r}(X, \beta)
\end{eqnarray}
forgetting the $(r+1)$-th marked point. Let $\mathcal L_i$
be the cotangent line bundle on 
$\overline{\mathfrak M}_{g, r}(X, \beta)$
associated to the $i$-th marked point. Put 
\begin{eqnarray*}
\psi_i = c_1(\mathcal L_i).
\end{eqnarray*}
For $\gamma_1, \ldots, \gamma_r \in H^*(X; \Q)$ and 
nonnegative integers $k_1, \ldots, k_r$, define
\begin{eqnarray}  \label{def_gw}
  \langle \tau_{k_1}(\gamma_1) \cdots 
    \tau_{k_r}(\gamma_r) \rangle_{g, \beta} 
= \int_{[\overline{\mathfrak M}_{g, r}(X, \beta)]^{\rm vir}}
    \prod_{i=1}^r \psi_i^{k_i} \text{\rm ev}_i^*(\gamma_i).
\end{eqnarray}
Define the {\it reduced} Gromov-Witten potential of $X$ by
\begin{eqnarray}  \label{rgwp}
  {\bf F}_{\text{GW}}' \left (X; u, v| 
    \prod_{i=1}^r \tau_{k_i}(\gamma_i) \right ) 
= \sum_{\beta \ne 0} \sum_{g \ge 0}
    \left \langle \prod_{i=1}^r \tau_{k_i}(\gamma_i) 
    \right \rangle_{g, \beta} \,\, u^{2g-2}v^\beta
\end{eqnarray}
omitting the constant maps. For $\beta \ne 0$, 
the {\it reduced partition function} 
\begin{eqnarray*}
{\bf Z}_{\text{GW}}'\left (X; u| 
\prod_{i=1}^r \tau_{k_i}(\gamma_i) \right )_\beta
\end{eqnarray*}
of degree-$\beta$ Gromov-Witten invariants is defined by 
setting:
\begin{eqnarray}  \label{rpf_gw1}
  1 + \sum_{\beta \ne 0} {\bf Z}_{\text{GW}}'\left (X; u| 
    \prod_{i=1}^r \tau_{k_i}(\gamma_i) \right )_\beta 
    \, v^\beta
\,\, = \,\,
  \text{exp} \,\, {\bf F}_{\text{GW}}' \left (X; u, v| 
    \prod_{i=1}^r \tau_{k_i}(\gamma_i) \right ).
\end{eqnarray}

Alternatively, let $\overline{\mathfrak M}_{g, r}'(X, \beta)$ 
be the moduli space of stable maps from 
{\it possibly disconnected} curves $C$ of genus-$g$ 
with $r$ marked points and with no collapsed connected 
components. Here the genus of a possibly disconnected 
curve $C$ is
\begin{eqnarray*}  
1 - \chi(\mathcal O_C) = 1 -\ell + \sum_{i=1}^\ell g_{C_i}
\end{eqnarray*}
where $C_1, \ldots, C_\ell$ denote all the connected 
components of $C$. For $\gamma_1, \ldots, \gamma_r \in 
H^*(X; \Q)$ and $k_1, \ldots, k_r \ge 0$, define
the reduced Gromov-Witten invariant by
\begin{eqnarray}  \label{def_rgw}
  \langle \tau_{k_1}(\gamma_1) \cdots 
    \tau_{k_r}(\gamma_r) \rangle_{g, \beta}' 
= \int_{[\overline{\mathfrak M}_{g, r}'(X, \beta)]^{\rm vir}}
    \prod_{i=1}^r \psi_i^{k_i} \text{\rm ev}_i^*(\gamma_i).
\end{eqnarray}
Then the reduced partition function of degree-$\beta$
invariants is also given by
\begin{eqnarray}  \label{rpf_gw2}
{\bf Z}_{\text{GW}}'\left (X; u| 
\prod_{i=1}^r \tau_{k_i}(\gamma_i) \right )_\beta
= \sum_{g \in \mathbb Z}
  \left \langle \prod_{i=1}^r \tau_{k_i}(\gamma_i) 
    \right \rangle_{g, \beta}' \,\, u^{2g-2}.
\end{eqnarray}

When $\dim (X)=3$, the expected dimensions of 
$\overline{\mathfrak M}_{g, r}(X, \beta)$ 
and $\overline{\mathfrak M}_{g, r}'(X, \beta)$ are
\begin{eqnarray}  \label{dim_gw}
-\int_{\beta} K_X + r.
\end{eqnarray}

\begin{remark}  \label{rmk1}
By the Fundamental Class Axiom, Divisor Axiom and
Dilation Axiom of the descendent Gromov-Witten invariants, 
if $\beta \ne 0$ and $\displaystyle{\int_{\beta} K_X} =0$, 
then
\begin{eqnarray*}
{\bf Z}_{\text{\rm GW}}'\left (X; u| 
\prod_{i=1}^r \tau_{k_i}(\gamma_i) \right )_\beta
\end{eqnarray*}
can be reduced to the case $r = 0$, i.e., 
to the reduced partition function
\begin{eqnarray}    \label{gw_r0}
{\bf Z}_{\text{\rm GW}}'\left (X; u \right )_\beta.
\end{eqnarray}
\end{remark}

\subsection{\bf The computations} 
\label{subsect_gw_comp}
\par
$\,$

We begin with the Gromov-Witten invariants of a smooth 
elliptic curve $E$. Let $d \ge 1$ and $[E] \in 
H_2(E; \mathbb Z)$ be the fundamental class. We use
\begin{eqnarray*}
\overline{\mathfrak M}_{g, r}(E, d), \;\;
\overline{\mathfrak M}_{g, r}'(E, d)
\end{eqnarray*}
to denote the moduli spaces 
$\overline{\mathfrak M}_{g, r}(E, d[E])$, 
$\overline{\mathfrak M}_{g, r}'(E, d[E])$ respectively.
The expected dimension of the moduli spaces 
$\overline{\mathfrak M}_{1, 0}(E, d)$
and $\overline{\mathfrak M}_{1, 0}'(E, d)$ is zero. So
\begin{eqnarray}  \label{def_E_gw}
\langle \rangle_{1, d[E]} 
&=&\text{deg} \big [\overline{\mathfrak M}_{1, 0}(E, d) 
   \big ]^{\text{vir}}, \\
\langle \rangle_{1, d[E]}' 
&=& \text{deg} \big [\overline{\mathfrak M}_{1, 0}'(E, d) 
   \big ]^{\text{vir}}.
\end{eqnarray}
Note that if $C$ is the (possibly disconnected) domain curve 
of a stable map in $\overline{\mathfrak M}_{1, 0}'(E, d)$,
then every connected component of $C$ must be of genus-$1$.
Therefore, as in (\ref{rgwp}), (\ref{rpf_gw1}) and 
(\ref{rpf_gw2}), we obtain the following relation:
\begin{eqnarray}  \label{rpf_E_gw}
 1 + \sum_{d=1}^{+ \infty} \langle \rangle_{1, d[E]}' \; v^d
\,\, = \,\,
 \text{exp} \; \sum_{d=1}^{+ \infty} 
 \langle \rangle_{1, d[E]} \; v^d.
\end{eqnarray}
By the Theorem~5 in \cite{OP1} (replacing $n$ and $q$ there
by $0$ and $v$ respectively),
\begin{eqnarray}  \label{thm5}
1 + \sum_{d=1}^{+ \infty} \langle \rangle_{1, d[E]}' \; v^d
\; = \; \frac{1}{\prod_{m=0}^{+\infty}(1-v^m)}.
\end{eqnarray}

In the rest of this section, we adopt the following 
notation.

\begin{notation}  \label{nota}
(i) Let $X = E \times S$ where $E$ is an elliptic curve and 
$S$ is a smooth surface. Let $\beta_0 \in H_2(X; \Z)$ be 
the fiber class of the fibration 
\begin{eqnarray*}
f: X = E \times S \to S.
\end{eqnarray*}
We use $K_X$ to denote both the canonical class and 
the canonical line bundle of $X$.

(ii) For $d \ge 0$, let $S^{[d]}$ be the Hilbert scheme 
which parametrizes the length-$d$ $0$-dimensional closed 
subschemes of the surface $S$. 

(iii) Fix $O \in E$ as the zero element for the group 
law on $E$. For $p \in E$, let
\begin{eqnarray}  \label{law}
\phi_p: E \to E
\end{eqnarray}
be the automorphism of $E$ defined via translation
$\phi_p(e) = p + e$. We have an action of $E$ on 
$X = E \times S$ via the automorphisms 
$\phi_p \times \id_S, \; p \in E$. 
\end{notation}

\begin{lemma} \label{lma_gw_d}
Let $X$ be from Notation~\ref{nota} and $d \ge 1$. 
Then, we have
\begin{eqnarray*}
\langle \rangle_{1, d \beta_0}' = \chi \big ( S^{[d]} \big ).
\end{eqnarray*}
\end{lemma}
\begin{proof}
First of all, let $\mathcal H_1^E$ be the rank-$1$ Hodge 
bundle over $\overline{\mathfrak M}_{1, 0}(E, d)$, i.e., 
\begin{eqnarray*}
\mathcal H_1^E = (\pi_{1, 0})_*\omega_{1, 0}
\end{eqnarray*}
where $\omega_{1, 0}$ is the relative dualizing sheaf of 
the forgetful map $\pi_{1, 0}$ in (\ref{pigr}).

Next, by the universal property of moduli spaces, we have
\begin{eqnarray} \label{lma_gw_d.1}
\overline{\mathfrak M}_{1, 0}(X, d \beta_0) \cong 
\overline{\mathfrak M}_{1, 0}(E, d) \times S.
\end{eqnarray}
By the definitions of virtual fundamental classes and 
the Hodge bundle,
\begin{eqnarray}   \label{lma_gw_d.2}
   \big [ \overline{\mathfrak M}_{1, 0}(X, d \beta_0) 
      \big ]^{\rm vir}
=  e \left ( \pi_1^* \big ( \mathcal H_1^E \big )^{\rm v} 
      \otimes \pi_2^*T_S \right ) 
   \, \cap \, \pi_1^*\big [ 
   \overline{\mathfrak M}_{1, 0}(E, d)\big ]^{\rm vir}
\end{eqnarray}
where $\pi_1$ and $\pi_2$ are the two projections of 
$\overline{\mathfrak M}_{1, 0}(X, d \beta_0)$ via 
the isomorphism (\ref{lma_gw_d.1}), and $e(\cdot)$ 
denotes the Euler class (or the top class). Note that
\begin{eqnarray*}  
  e \big ( \pi_1^* \big ( \mathcal H_1^E \big )^{\rm v}
      \otimes \pi_2^*T_S \big )
= \pi_2^*e(S) + \pi_2^*K_S \cdot 
      \pi_1^*c_1\big ( \mathcal H_1^E \big )
      + \pi_1^*c_1\big ( \mathcal H_1^E \big )^2.
\end{eqnarray*}
By (\ref{lma_gw_d.2}), $\langle \rangle_{1, d \beta_0}
= \chi(S) \cdot \langle \rangle_{1, d[E]}$. 
Therefore, we obtain
\begin{eqnarray}   \label{lma_gw_d.3}
   1 + \sum_{d=1}^{+\infty} \langle \rangle_{1, d \beta_0}' 
     \; v^d
&=&\text{\rm exp} \, \sum_{d=1}^{+\infty} 
     \langle \rangle_{1, d \beta_0} \; v^d   
     \nonumber  \\
&=&\text{\rm exp} \, \left ( \chi(S) \cdot \sum_{d=1}^{+\infty} 
     \langle \rangle_{1, d[E]} \; v^d \right )  
     \nonumber  \\
&=&\left ( 1 + \sum_{d=1}^{+\infty} 
     \langle \rangle_{1, d[E]}' \; v^d \right )^{\chi(S)}  
     \nonumber  \\
&=&\frac{1}{\prod_{m=0}^{+\infty}(1-v^m)^{\chi(S)}}
\end{eqnarray}
by (\ref{rpf_E_gw}) and (\ref{thm5}). By G\" ottsche's 
formula in \cite{Got} for $\chi \big ( S^{[d]} \big )$, we have
\begin{eqnarray*}
  \sum_{d=0}^{+\infty} \chi \big ( S^{[d]} \big ) \, v^d
= \frac{1}{\prod_{m=0}^{+\infty}(1-v^m)^{\chi(S)}}.
\end{eqnarray*}
Combining this with (\ref{lma_gw_d.3}), we conclude that 
$\langle \rangle_{1, d \beta_0}' = \chi \big ( S^{[d]} \big )$.
\end{proof}

Let $X$ be from Notation~\ref{nota} and 
$\beta \in H_2(X; \Z) \backslash \{0\}$.
For any $p \in E$,
\begin{eqnarray}  \label{fix_beta}
(\phi_p \times \text{\rm Id}_S)_*\beta = \beta
\end{eqnarray}
since $\{\phi_p \times \text{\rm Id}_S\}_{p \in E}$
form a connected algebraic family of automorphims of $X$.
Thus the algebraic group $E$ acts on the stack
of $r$-pointed degree-$\beta$ stable maps to $X$ 
(see \cite{Kon}). The universal properties of moduli spaces 
imply that there is a corresponding action of $E$ on 
the moduli space $\overline{\mathfrak M}_{g, r}(X,\beta)$.
For $p \in E$, let
\begin{eqnarray*}
\Psi_p: \overline{\mathfrak M}_{g, r}(X,\beta) \to 
        \overline{\mathfrak M}_{g, r}(X,\beta)
\end{eqnarray*}
be the corresponding automorphism. Then we see that 
the automorphism $\Psi_p$ maps a point 
$[\mu: (C; w_1, \ldots, w_r) \to X] \in 
\overline{\mathfrak M}_{g, r}(X,\beta)$ to the point
\begin{eqnarray}   \label{Psip}
[(\phi_p \times \text{\rm Id}_S) \circ \mu: 
(C; w_1, \ldots, w_r) \to X] 
\in \overline{\mathfrak M}_{g, r}(X,\beta).
\end{eqnarray}

\begin{lemma} \label{lma_eaction}
With the notation as above, the algebraic group $E$ acts 
without fixed points on 
$\overline{\mathfrak M}_{g, r}(X,\beta)$
if $\beta \neq d \beta_0$, or $r \ge 1$, or $g \neq 1$.
\end{lemma}
\begin{proof}
Assume that $[\mu: (C; w_1, \ldots, w_r) \to X] 
\in \overline{\mathfrak M}_{g, r}(X,\beta)$ is fixed by 
the action of $E$. By definition, for every $p \in E$,
there is an automorphism $\tau_p$ of $C$ such that 
\begin{eqnarray}  \label{lma_eaction.1}
\mu \circ \tau_p = (\phi_p \times \text{\rm Id}_S) \circ \mu
\end{eqnarray}
and $\tau_p(w_i) = w_i$ for all $1 \le i \le r$.
In particular, for every $p \in E$, we have
\begin{eqnarray*}
\mu(C) = (\phi_p \times \text{\rm Id}_S)
\big ( \mu(C) \big ).
\end{eqnarray*}
So $\mu(C)$ is a fiber of the elliptic fibration $f$,
and $\beta = d \beta_0$ for some $d \ge 1$.
By our assumption, either $r \ge 1$ or $g \ge 2$.
By (\ref{lma_eaction.1}), we get
\begin{eqnarray}  \label{lma_eaction.2}
\mu \circ \tau_p (C) = \phi_p \big ( \mu(C) \big ).
\end{eqnarray}
Since $\phi_p$ acts freely on the fiber $\mu(C)$, 
(\ref{lma_eaction.2}) implies that the automorphisms 
$\tau_p$ of the marked curve $(C; w_1, \ldots, w_r)$ 
are different for different points $p \in E$. 
Hence the automorphism group of the marked curve 
$(C; w_1, \ldots, w_r)$ is infinite. This is impossible
since either $g \ge 2$ or $g = 1$ and $r \ge 1$.
\end{proof}

\begin{lemma} \label{lma_gw=0}
Let $\beta \in H_2(X; \Z) \backslash \{0\}$. 
Assume that $\gamma_1, \ldots, \gamma_r \in f^*H^*(S; \Q)
\subset H^*(X; \Q)$. If $\beta \neq d \beta_0$, 
or $r \ge 1$, or $g \neq 1$, then we have
\begin{eqnarray*}  
\langle \tau_{k_1}(\gamma_1) \cdots 
\tau_{k_r}(\gamma_r) \rangle_{g, \beta}' = 0.
\end{eqnarray*}
\end{lemma}
\begin{proof}
First of all, note that it suffices to show that
\begin{eqnarray}  \label{lma_gw=0.1}
\langle \tau_{k_1}(\gamma_1) \cdots 
\tau_{k_r}(\gamma_r) \rangle_{g, \beta} = 0
\end{eqnarray}
if $\beta \neq d \beta_0$, or $r \ge 1$, or $g \neq 1$. 
In the following, we prove (\ref{lma_gw=0.1}).

By Lemma~\ref{lma_eaction}, $E$ acts without fixed 
points on $\overline{\mathfrak M}_{g, r}(X,\beta)$. 
Since $E$ is an elliptic curve, any proper algebraic 
subgroup is finite. Thus the stabilizer of any point for 
the $E$ action on $\overline{\mathfrak M}_{g, r}(X,\beta)$
is finite. Since $\overline{\mathfrak M}_{g, r}(X,\beta)$ 
is finite type, the order of the stabilizer subgroup 
at any point is bounded by some number $N$. Thus, 
if $G$ is a cyclic subgroup of $E$ of prime order 
$p > N$, then $G$ acts freely on 
$\overline{\mathfrak M}_{g, r}(X,\beta)$. We fix such 
a cyclic subgroup $G$ of $E$ in the rest of the proof.

The complex $R(\pi_{g,r})_* ({\rm ev}_{r+1})^*T_X$ from 
(\ref{pot_gw}) is equivariant for the action of 
any algebraic automorphism group of $X$.
Thus for some positive integer $m$ (independent of $G$), 
the cycle $m \big [\overline{\mathfrak M}_{g,
r}(X,\beta) \big ]^{\rm vir}$
defines an element of the integral equivariant Borel-Moore 
homology group $H_*^G \big (
\overline{\mathfrak M}_{g, r}(X, \beta) \big )$.
Likewise if $\gamma_i \in f^*H^*(S; \Q)$, 
then the cycle $\gamma_i$ is invariant under 
the action of $E$ on $X$. Hence some positive multiple 
$m_i \gamma_i$ defines an element of $H^*_G(X)$, 
where $m_i$ is independent of $G$. Note from (\ref{Psip}) 
that the evaluation map
${\rm ev}_i: \overline{\mathfrak M}_{g,r}(X,\beta) \to X$
is $G$-equivariant, so the pullback ${\rm ev}_i^*(m_i \gamma_i)$
determines an element of 
$H^*_G(\overline{\mathfrak M}_{g,r}(X,\beta))$.
In addition, the cotangent line bundles ${\mathcal L}_i$ 
($1 \leq i \leq r$) over 
$\overline{\mathfrak M}_{g,r}(X,\beta)$
are equivariant for the action of $G$. It follows 
from the definition (\ref{def_gw}) that the cycle 
\begin{eqnarray*}
m m_1 \cdots m_r \; \langle \tau_{k_1}(\gamma_1) 
\cdots \tau_{k_r}(\gamma_r) \rangle_{g, \beta}
\end{eqnarray*}
defines an element in the degree-$0$ Borel-Moore
homology $H^G_0 \big (
\overline{\mathfrak M}_{g, r}(X,\beta) \big )$. 

Since $G$ is a cyclic subgroup of order $p$ which acts 
freely on $\overline{\mathfrak M}_{g, r}(X,\beta)$,
any element of $H^G_0 \big (\overline{\mathfrak
M}_{g, r}(X,\beta) \big )$ is represented by 
a $G$-invariant $0$-cycle whose degree is a multiple 
of $p$ (possibly 0). Since $p$ can be taken to be 
arbitrarily large,
\begin{eqnarray*}
m m_1 \cdots m_r \; \langle \tau_{k_1}(\gamma_1) 
\cdots \tau_{k_r}(\gamma_r) \rangle_{g, \beta}  = 0.
\end{eqnarray*}
Therefore, $\langle \tau_{k_1}(\gamma_1) 
\cdots \tau_{k_r}(\gamma_r) \rangle_{g, \beta} = 0$.
This completes the proof of (\ref{lma_gw=0.1}).
\end{proof}

We define the cohomology degree $|\gamma| = \ell$ when 
$\gamma \in H^{\ell}(X; \Q)$.

\begin{proposition} \label{prop_gw2}
Let $\beta \in H_2(X; \Z) \backslash \{0\}$. Assume 
$\displaystyle{\int_{\beta} K_X} 
= \displaystyle{\int_{\beta} f^*K_S} =0$. Then,
\begin{eqnarray*}
&&  {\bf Z}_{\text{\rm GW}}'\left (X; u| 
   \prod_{i=1}^r \tau_{0}(\gamma_i) \right )_\beta \\
&=&\left \{
   \begin{array}{ll}
     \prod_{i=1}^r \int_{\beta} \gamma_i 
       \cdot \chi \big ( S^{[d]} \big )
       &\text{if} \; |\gamma_i| = 2 \; 
        \text{for every} \; i \text{ and } \beta = d\beta_0
     \text{ for some } d \ge 1;\\
     0 &\text{otherwise}.
   \end{array}
  \right.
\end{eqnarray*}
\end{proposition}
\begin{proof}
By (\ref{dim_gw}) and the degree condition on 
Gromov-Witten invariants, 
\begin{eqnarray*}
\sum_{i=1}^r |\gamma_i| = 2r.
\end{eqnarray*}
By the Fundamental Class Axiom and Divisor Axiom of 
Gromov-Witten invariants,
\begin{eqnarray*}  
\langle \tau_{0}(\gamma_1) \cdots 
\tau_{0}(\gamma_r) \rangle_{g, \beta} 
= \left \{
   \begin{array}{ll}
     \prod_{i=1}^r \int_{\beta} \gamma_i
       \cdot \langle  \rangle_{g, \beta}
       &\text{\rm if} \; |\gamma_i| = 2 \; 
        \text{\rm for every} \; i;\\
     0 &\text{\rm otherwise}.
   \end{array}
  \right.
\end{eqnarray*}
So by Lemma~\ref{lma_gw_d} and by taking $r=0$ in 
(\ref{lma_gw=0.1}), we conclude that
\begin{eqnarray*}  
&&\langle \tau_{0}(\gamma_1) \cdots 
  \tau_{0}(\gamma_r) \rangle_{g, \beta}'  \\
&=&\left \{
   \begin{array}{ll}
     \prod_{i=1}^r \int_{\beta} \gamma_i
       \cdot \chi \big ( S^{[d]} \big )
       &\text{\rm if} \; |\gamma_i| = 2 \; 
        \text{\rm for every} \; i, g=1, \beta = d\beta_0;\\
     0 &\text{\rm otherwise}.
   \end{array}
  \right.
\end{eqnarray*}
Now our proposition follows directly from the identity 
(\ref{rpf_gw2}).
\end{proof}

\begin{proposition} \label{prop_gw1}
Let $X$ be from Notation~\ref{nota} and 
$\beta \in H_2(X; \Z) \backslash \{0\}$.
Assume that $\gamma_1, \ldots, \gamma_r \in f^*H^*(S; \Q)
\subset H^*(X; \Q)$. Then,
\begin{eqnarray*}
  {\bf Z}_{\text{\rm GW}}'\left (X; u| 
   \prod_{i=1}^r \tau_{k_i}(\gamma_i) \right )_\beta
= \left \{
   \begin{array}{ll}
     \chi \big ( S^{[d]} \big ) &if \; r = 0 \text{ and }
       \beta = d\beta_0 \text{ with } d \ge 1;\\
     0 &otherwise.
   \end{array}
  \right.
\end{eqnarray*}
\end{proposition}
\begin{proof}
Follows from the identity (\ref{rpf_gw2}), 
Lemma~\ref{lma_gw_d} and Lemma~\ref{lma_gw=0}.
\end{proof}

\section{\bf Donaldson-Thomas theory} 
\label{sect_dt}
\subsection{\bf Donaldson-Thomas invariants} 
\label{subsect_dt_inv}
\par
$\,$

Let $X$ be a smooth projective complex threefold. 
For a fixed class $\beta \in H_2(X; \Z)$
and a fixed integer $n$, following the definition and
notation in \cite{MNOP1, MNOP2},
we define $\mathfrak I_n(X, \beta)$ to be the moduli space 
parametrizing the ideal sheaves $I_Z$ of 
$1$-dimensional closed subschemes $Z$ of $X$ satisfying
the conditions:
\begin{eqnarray}   \label{In}
\chi(\mathcal O_Z)=n, \qquad [Z] = \beta
\end{eqnarray}
where $[Z]$ is the class associated to the dimension-$1$ 
component (weighted by their intrinsic multiplicities) of $Z$.
Note that $\mathfrak I_n(X, \beta)$ is a special case of
the moduli spaces of Gieseker semistable torsion-free sheaves
over $X$. 
When the anti-canonical divisor $-K_X$ 
is effective, perfect obstruction theories on the moduli spaces
$\mathfrak I_n(X, \beta)$ have been constructed in \cite{Tho}.
This result has been generalized in \cite{MP}.
By the Lemma~1 in \cite{MNOP2}, the virtual dimension of 
$\mathfrak I_n(X, \beta)$ is
\begin{eqnarray}   \label{dim_In}
- \int_{\beta} K_X.
\end{eqnarray}

The Donaldson-Thomas invariant is defined via integration
against the virtual fundamental class 
$[\mathfrak I_n(X, \beta)]^{\text{vir}}$ of the moduli space
$\mathfrak I_n(X, \beta)$. More precisely,
let $\gamma \in H^\ell(X; \Q)$ and $\mathcal I$ be 
the universal ideal sheaf over 
$\mathfrak I_n(X, \beta) \times X$. Let
\begin{eqnarray}  \label{def_opr1}
\och_{k+2}(\gamma): \; 
H_*\big ( \mathfrak I_n(X, \beta); \Q \big ) \to
H_{*-2k+2-\ell} \big ( \mathfrak I_n(X, \beta); \Q \big )
\end{eqnarray}
be the operation on the homology of $\mathfrak I_n(X, \beta)$
defined by 
\begin{eqnarray}  \label{def_opr2}
\och_{k+2}(\gamma)(\xi) = \pi_{1*} \big ( \ch_{k+2}(\mathcal I) 
\cdot \pi_2^*\gamma \cap \pi_1^*\xi \big )
\end{eqnarray}
where $\pi_1$ and $\pi_2$ be the two projections on 
$\mathfrak I_n(X, \beta) \times X$. Define
\begin{eqnarray}  \label{dt_desc}
& &\langle \wtau_{k_1}(\gamma_1) \cdots 
   \wtau_{k_r}(\gamma_r) \rangle_{n, \beta}  \nonumber   \\
&=&\int_{[\mathfrak I_n(X, \beta)]^{\text{vir}}}
   \prod_{i=1}^r (-1)^{k_i+1}\och_{k_i+2}(\gamma_i)  
   \nonumber   \\
&=&(-1)^{k_1+1}\och_{k_1+2}(\gamma_1) \circ \cdots \circ
  (-1)^{k_r+1}\och_{k_r+2}(\gamma_r)
  \big ( [\mathfrak I_n(X, \beta)]^{\text{vir}} \big ).
\end{eqnarray}
The partition function for these descendent Donaldson-Thomas 
invariants is
\begin{eqnarray}  \label{p_func}
{\bf Z}_{\text{DT}} \left (X; q \; | \prod_{i=1}^r
\wtau_{k_i}(\gamma_i) \right )_\beta = 
\sum_{n \in \Z} \left \langle \wtau_{k_1}(\gamma_1) \cdots 
\wtau_{k_r}(\gamma_r) \right \rangle_{n, \beta} \; q^n.
\end{eqnarray}

The partition function for the degree-$0$ Donaldson-Thomas 
invariants of $X$ is
\begin{eqnarray}  \label{deg0}
{\bf Z}_{\text{DT}}(X; q)_0 = M(-q)^{\chi(X)}
\end{eqnarray}
by \cite{JLi, BF2} (this formula was conjectured in 
\cite{MNOP1, MNOP2}), where 
\begin{eqnarray*}
M(q) = \prod_{n=1}^{+\infty} \frac{1}{(1-q^n)^n}
\end{eqnarray*}
is the McMahon function. The {\it reduced partition function}
is defined to be
\begin{eqnarray}  \label{rpf_dt}
   {\bf Z}_{\text{DT}}' \left (X; q \; | \prod_{i=1}^r
     \wtau_{k_i}(\gamma_i) \right )_\beta
&=&\frac{{\bf Z}_{\text{DT}} \left (X; q \; | \prod_{i=1}^r
     \wtau_{k_i}(\gamma_i) \right )_\beta}{
     {\bf Z}_{\text{DT}}(X, q)_0}         \nonumber   \\
&=&\frac{{\bf Z}_{\text{DT}} \left (X; q \; | \prod_{i=1}^r
     \wtau_{k_i}(\gamma_i) \right )_\beta}{M(-q)^{\chi(X)}}.   
\end{eqnarray}

In the next two lemmas, we study the operators
$\och_{2}(\gamma)$ and $\och_{3}(1_X)$ respectively, 
where $1_X \in H^*(X; \Q)$ is the fundamental cohomology 
class. The results will be used in 
Subsect.~\ref{subsect_dt_comp}. 
Note that the first lemma is the analogue to 
the Fundamental Class Axiom and Divisor Axiom of 
Gromov-Witten invariants, while the second one is 
the analogue to the Dilaton Axiom of 
Gromov-Witten invariants. By (\ref{def_opr1}),
\begin{eqnarray*}  
\och_{2}(\gamma): \; 
  H_b \big ( \mathfrak I_n(X, \beta); \Q \big ) \to
  H_{b-2+|\gamma|} \big ( \mathfrak I_n(X, \beta); \Q \big ), 
\end{eqnarray*}
\begin{eqnarray*}  
\och_{3}(1_X): \; 
  H_b \big ( \mathfrak I_n(X, \beta); \Q \big ) \to
  H_{b} \big ( \mathfrak I_n(X, \beta); \Q \big ).
\end{eqnarray*}
Let $cl: A_* \big ( \mathfrak I_n(X, \beta) \big ) \otimes \Q
\to H_* \big ( \mathfrak I_n(X, \beta); \Q \big )$ be 
the cycle map. Put
\begin{eqnarray*}  
H_*^{\rm alg} \big ( 
\mathfrak I_n(X, \beta) \big ) = \text{\rm im} (cl). 
\end{eqnarray*}

\begin{lemma} \label{dt_fund_div}
{\rm (i)} Let $\beta \in H_2(X; \Z)$ and 
$\gamma \in H^\ell(X; \Q)$. Then,
\begin{eqnarray*}
  \och_{2}(\gamma)|_{H_*^{\rm alg}(\mathfrak I_n(X, \beta))} 
= \left \{
   \begin{array}{ll}
     0 &\text{ if } \ell = 0 \text{ or } 1;\\
     -\int_{\beta} \gamma \cdot \id &\text{ if } \ell =2.
   \end{array}
  \right.
\end{eqnarray*}

{\rm (ii)} If the moduli space $\mathfrak I_n(X, \beta)$
is smooth, then
\begin{eqnarray*}
  \och_{2}(\gamma) 
= \left \{
   \begin{array}{ll}
     0 &\text{ if } \ell = 0 \text{ or } 1;\\
     -\int_{\beta} \gamma \cdot \id &\text{ if } \ell =2.
   \end{array}
  \right.
\end{eqnarray*}
\end{lemma}
\begin{proof}
(i) Let $\mathfrak I = \mathfrak I_n(X,\beta)$. 
By \cite{FG}, there is a proper morphism 
\begin{eqnarray*}
p: \; \W{\mathfrak I} \to \mathfrak I
\end{eqnarray*}
with $\W{\mathfrak I}$ smooth and 
$p_*\colon H_*^{\rm alg}(\W{\mathfrak I}) \to 
H_*^{\rm alg}(\mathfrak I)$ surjective. 
Such a morphism $p$ is called a {\it nonsingular envelope} 
(see p.299 of \cite{FG}).
Let $\tilde{\pi}_1$ and $\tilde{\pi}_2$ be the projections
from $\W{\mathfrak I} \times X$ to the first and 
second factors respectively.

Let $\xi \in H_*^{\rm alg}(\mathfrak I)$. Then
$\xi = p_* \tilde{\xi}$ for some 
$\w \xi \in H_*^{\rm alg}(\W{\mathfrak I})$. Define
\begin{eqnarray}   \label{dan.1}
\tilde{\och}_2(\gamma)(\tilde{\xi}) = \tilde{\pi}_{1*}
\left ( \ch_2 \big ((p \times \id_X)^* \mathcal I \big ) \, 
\tilde{\pi}_2^* \gamma \cap \tilde{\pi}_1^* \tilde{\xi} 
\right )
\end{eqnarray}
where $\mathcal I$ denotes the universal ideal sheaf over 
$\mathfrak I \times X$. Using the projection formula and 
the fact that $(p \times \id_X)_*\tilde{\pi}_1^* \tilde{\xi} 
= {\pi}_1^* p_*\tilde{\xi} = {\pi}_1^* \xi$, we have
\begin{eqnarray}  \label{dan.2}
   p_* \big (\tilde{\och}_2(\gamma)(\tilde{\xi}) \big )
&=&p_*\tilde{\pi}_{1*}
   \left ( \ch_2 \big ((p \times \id_X)^* \mathcal I \big ) \, 
   \tilde{\pi}_2^* \gamma \cap \tilde{\pi}_1^* \tilde{\xi} 
   \right )          \nonumber   \\
&=&\pi_{1*}(p \times \id_X)_*
   \left ( (p \times \id_X)^* \big ( \ch_2(\mathcal I) \, 
   \pi_2^* \gamma  \big ) \cap \tilde{\pi}_1^* \tilde{\xi} 
   \right )          \nonumber    \\
&=&\pi_{1*}\left ( \ch_2(\mathcal I) \, 
   \pi_2^* \gamma \cap (p \times \id_X)_*\tilde{\pi}_1^*
   \tilde{\xi} \right )         \nonumber     \\
&=&\och_2(\gamma)(\xi).
\end{eqnarray}
Since $\W{\mathfrak I}$ is smooth, the Poincar\'e duality
holds and we see from (\ref{dan.1}) that
\begin{eqnarray*}
\tilde{\och}_2(\gamma)(\tilde{\xi}) 
= \tilde{\pi}_{1*} \big (\ch_2((p \times \id_X)^* \mathcal I)
\tilde{\pi}_2^*\gamma \big ) \cap \w \xi
\end{eqnarray*}
where $\tilde{\pi}_{1*} \big (\ch_2((p \times \id_X)^*
\mathcal I) \tilde{\pi}_2^*\gamma \big )$ is 
the cohomology class Poincar\'e dual to 
\begin{eqnarray*}
\tilde{\pi}_{1*} \left (\ch_2((p \times \id_X)^* \mathcal I)
\tilde{\pi}_2^*\gamma \cap [\W{\mathfrak I} \times X] \right).
\end{eqnarray*}
Thus by (\ref{dan.2}), to prove the lemma, it suffices to 
show that
\begin{eqnarray}   \label{dan.3}
  \tilde{\pi}_{1*}\left(\ch_2((p \times \id_X)^*\mathcal I) 
  \tilde{\pi}_2^*\gamma \cap [\W{\mathfrak I} \times X]\right)
= \left \{
   \begin{array}{ll}
     0 &\text{ if } \ell = 0 \text{ or } 1;\\
     -\int_{\beta} \gamma \cdot [\W{\mathfrak I}] 
           &\text{ if } \ell =2.
   \end{array}
  \right.
\end{eqnarray}

Let $\mathcal Z \subset \mathfrak I \times X$
be the universal closed subscheme. Set-theoretically, 
\begin{eqnarray*}  
\mathcal Z = \{ (I_Z, x) \in \mathfrak I 
\times X|\; x \in \Supp(Z) \}.
\end{eqnarray*}
Let $\W{\mathcal Z} = (p \times \id_X)^{-1}\mathcal Z$.
Then, $\mathcal I = I_{\mathcal Z}$,
$(p \times \id_X)^*\mathcal I = (p \times \id_X)^*
I_{\mathcal Z} = I_{\W{\mathcal Z}}$, and 
\begin{eqnarray}    \label{dan.4}
\ch_2((p \times \id_X)^* \mathcal I)
= \ch_2(I_{\W{\mathcal Z}}) = -c_2(I_{\W{\mathcal Z}})
= c_2(\mathcal O_{\W{\mathcal Z}}).
\end{eqnarray}

If $\beta = 0$, then $\mathcal Z$ is of codimension-$3$ in
$\mathfrak I \times X$, and $\W{\mathcal Z}$ is of 
codimension-$3$ in $\W{\mathfrak I} \times X$ as well. 
By (\ref{dan.4}), $\ch_2((p \times \id_X)^* \mathcal I)
= 0$. Therefore, (\ref{dan.3}) holds.

Next, we assume $\beta \neq 0$. Then, $\mathcal Z$ is of 
codimension-$2$ in $\mathfrak I \times X$, 
and $\W{\mathcal Z}$ is of codimension-$2$ in 
$\W{\mathfrak I} \times X$. By (\ref{dan.4}), 
$\ch_2((p \times \id_X)^* \mathcal I)
= -[\W{\mathcal Z}]$. So
\begin{eqnarray}   \label{dan.5}
\tilde{\pi}_{1*}
\left(\ch_2((p \times \id_X)^* \mathcal I) 
\tilde{\pi}_2^*\gamma \cap [\W{\mathfrak I} \times X]\right) 
= -\tilde{\pi}_{1*} \big ([\W{\mathcal Z}] \cdot
  \tilde{\pi}_2^*\gamma \big ).
\end{eqnarray}
When $\ell = 0$ or $1$, we get $\tilde{\pi}_{1*} \big 
([\W{\mathcal Z}] \cdot \tilde{\pi}_2^*\gamma \big ) = 0$
by degree reason. Hence (\ref{dan.3}) holds.

We are left with the case $\ell = 2$. In this case, 
$\tilde{\pi}_{1*} \big ([\W{\mathcal Z}] \cdot
\tilde{\pi}_2^*\gamma \big )$ is a multiple of $[\W{\mathfrak I}]$.
Let $m$ be the multiplicity, and $\w w \in \W{\mathfrak I}$ 
be a point. Then, we have
\begin{eqnarray*}
m = \deg \; \big ( [\W{\mathcal Z}] \cdot 
\tilde{\pi}_2^*\gamma \big ) |_{\{\w w\} \times X}
= \int_{\beta} \gamma.
\end{eqnarray*}
Therefore, we conclude from (\ref{dan.5}) that 
(\ref{dan.3}) holds when $\ell = 2$.

(ii) Follows from the proof of (i) by taking $\W{\mathfrak I}
= {\mathfrak I}$ and $p = \id_{\mathfrak I}$.
\end{proof}

\begin{lemma}   \label{ch3}
{\rm (i)} Let $\beta \in H_2(X; \Z)$. Then, we have
\begin{eqnarray*}
  \och_3(1_X)|_{H_*^{\rm alg}(\mathfrak I_n(X, \beta))} 
= - \left (n + \int_\beta K_X \right ) \cdot \id.
\end{eqnarray*}

{\rm (ii)} If the moduli space $\mathfrak I_n(X, \beta)$
is smooth, then
\begin{eqnarray}   \label{ch3.1}
\och_3(1_X) = - \left (n + \int_\beta K_X \right ) \cdot \id.
\end{eqnarray}
\end{lemma}
\begin{proof}
Note that (i) follows from the proof of (ii) and 
the similar trick of using a nonsingular envelope as 
in the proof of Lemma~\ref{dt_fund_div}~(i).
To prove (ii), we adopt the notation in (\ref{def_opr2}). 
Using the projection formula, we get
\begin{eqnarray}   \label{ch3.2}
\och_3(1_X)(\xi) = \pi_{1*} \big ( \ch_{3}(\mathcal I) 
\cap \pi_1^*\xi \big ) = \pi_{1*} \ch_{3}(\mathcal I) 
\cdot \xi
\end{eqnarray}
since our moduli space $\mathfrak I_n(X, \beta)$ is smooth. 
Note that $\pi_{1*} \ch_{3}(\mathcal I)$ is a multiple of 
the fundamental cycle of $\mathfrak I_n(X, \beta)$.
Let $m$ be the multiplicity. Then,
\begin{eqnarray*}
m = \deg \; \ch_{3}(\mathcal I)|_{[I_Z] \times X}
  = \deg \; \ch_{3}(I_Z)
= - \deg \; \ch_3(\mathcal O_Z)
= - \frac{1}{2} \deg \; c_3(\mathcal O_Z)
\end{eqnarray*}
where $[I_Z]$ denotes a point in $\mathfrak I_n(X, \beta)$.
Since $c_1(\mathcal O_Z) = 0$ and $c_2(\mathcal O_Z) = -[Z]
=-\beta$, we see from (\ref{In}) and 
the Hirzebruch-Riemann-Roch Theorem that
\begin{eqnarray*}
  m
= - \frac{1}{2} \deg \; c_3(\mathcal O_Z)
= - \left (n + \int_\beta K_X \right ).
\end{eqnarray*}
Now combining this with (\ref{ch3.2}), we immediately 
obtain formula (\ref{ch3.1}).
\end{proof}

\begin{remark} \label{rmk_och2}
Let $\beta \in H_2(X; \Z)$ and $\gamma \in H^\ell(X; \Q)$. 
We expect that both Lemma~\ref{dt_fund_div} and Lemma~\ref{ch3} 
can be sharpened, i.e., we expect in general that
\begin{eqnarray*}
  \och_{2}(\gamma)
&=& \left \{
   \begin{array}{ll}
     0 &\text{ if } \ell = 0 \text{ or } 1;\\
     -\int_{\beta} \gamma \cdot \id &\text{ if } \ell =2;
   \end{array}
  \right.  \\  
\och_3(1_X) &=& - \left (n + \int_\beta K_X \right ) \cdot \id.
\end{eqnarray*}
\end{remark}

\subsection{\bf The computations} 
\label{subsect_dt_comp}
\par
$\,$

In the rest of this section, we adopt the notation in
Notation~\ref{nota}. We begin with the case when 
$n=0$ and $\beta = d \beta_0$ with $d \ge 0$. Note that
\begin{eqnarray}   \label{iso_hil}
\mathfrak I_0(X, d \beta_0) \cong S^{[d]}.
\end{eqnarray}
However, the expected dimension of 
$\mathfrak I_0(X, d \beta_0)$ is zero by (\ref{dim_In}).

\begin{lemma} \label{obstr}
{\rm (i)} The obstruction bundle over the moduli space 
$\mathfrak I_0(X, d \beta_0) \cong S^{[d]}$ is isomorphic 
to the tangent bundle $T_{S^{[d]}}$ of 
the Hilbert scheme $S^{[d]}$.

{\rm (ii)} The Donaldson-Thomas invariant 
$\langle \rangle_{0, d \beta_0}$ is equal to 
$\chi \big (S^{[d]} \big )$.
\end{lemma}
\begin{proof}
It is clear that (ii) follows from (i). To prove (i), let 
\begin{eqnarray*}
\psi = \id_{S^{[d]}} \times f: S^{[d]} \times X 
\to S^{[d]} \times S
\end{eqnarray*}
and $\phi: S^{[d]} \times S \to 
S^{[d]}$ be the projections. Let $\pi = \phi \circ \psi: 
S^{[d]} \times X \to S^{[d]}$. Let $\mathcal J$ be 
the universal ideal sheaf over $S^{[d]} \times S$.
Then the universal ideal sheaf over 
\begin{eqnarray*}
\mathfrak I_0(X, d \beta_0) \times X 
\cong S^{[d]} \times X
\end{eqnarray*}
is $\mathcal I = \psi^*\mathcal J$. The Zariski tangent
bundle and obstruction bundle over the moduli space
$\mathfrak I_0(X, d \beta_0) \cong S^{[d]}$ are given by 
the rank-$2d$ bundles
\begin{eqnarray*}
\mathcal Ext^1_\pi(\psi^*\mathcal J, \psi^*\mathcal J)_0, \;
\mathcal Ext^2_\pi(\psi^*\mathcal J, \psi^*\mathcal J)_0
\end{eqnarray*}
respectively (see, for instance, the Theorem~3.28 in \cite{Tho}
for the obstruction bundle). Here $\mathcal Ext^*_\pi$
denotes the right derived functors of $\mathcal Hom_\pi 
= \pi_* \mathcal Hom$. We claim
\begin{eqnarray}   \label{local.1}
\mathcal Ext^1_\pi(\psi^*\mathcal J, \psi^*\mathcal J)_0
\cong
\mathcal Ext^2_\pi(\psi^*\mathcal J, \psi^*\mathcal J)_0.
\end{eqnarray}

In the following, we will prove the local version of 
(\ref{local.1}), i.e., for every point $I_{f^*\xi} \in 
\mathfrak I_0(X, d \beta_0)$ with $\xi \in S^{[d]}$,
we show that there exists a canonical isomorphism:
\begin{eqnarray} \label{obstr.1}
Ext^1(I_{f^*\xi}, I_{f^*\xi})_0 \cong
Ext^2(I_{f^*\xi}, I_{f^*\xi})_0.
\end{eqnarray}
The argument for the global version (\ref{local.1}) follows 
from that for the local version (\ref{obstr.1}) and 
the isomorphisms via relative duality 
(see the Proposition~8.14 in \cite{LeP}):
\begin{eqnarray*}   
       \mathcal Ext^2_\pi(\psi^*\mathcal J, \psi^*\mathcal J)_0
&\cong&\mathcal Ext^1_\pi(\psi^*\mathcal J, \psi^*\mathcal J
       \otimes \w \rho^*K_S)_0^{\rm v}            \\
       \mathcal Ext^1_\phi(\mathcal J, \mathcal J)_0
&\cong&\mathcal Ext^1_\phi(\mathcal J, \mathcal J
       \otimes \rho^*K_S)_0^{\rm v}              
\end{eqnarray*}
where $\w \rho: S^{[d]} \times X = S^{[d]} \times S \times E
\to S$ and $\rho: S^{[d]} \times S \to S$ are the projections.

Here is an outline for (\ref{obstr.1}). 
We apply the Serre duality twice: once on $X$ with 
\begin{eqnarray*} 
       Ext^2(I_{f^*\xi}, I_{f^*\xi})_0 
&\cong&Ext^1(I_{f^*\xi}, I_{f^*\xi} \otimes K_X)_0^{\rm v} \\
&\cong&Ext^1(I_{f^*\xi}, I_{f^*\xi} \otimes f^*K_S)_0^{\rm v},
\end{eqnarray*}
and the other on $S$ with $Ext^1(I_{\xi}, I_{\xi})_0 \cong
Ext^1(I_{\xi}, I_{\xi} \otimes K_S)_0^{\rm v}$. 
Note from (\ref{iso_hil}) that
\begin{eqnarray} \label{obstr.2}
Ext^1(I_{f^*\xi}, I_{f^*\xi})_0 \cong
Ext^1(I_{\xi}, I_{\xi})_0.
\end{eqnarray}
The main part of our argument is to prove that there is 
a natural isomorphism:
\begin{eqnarray*} 
Ext^1(I_{f^*\xi}, I_{f^*\xi} \otimes f^*K_S)_0 \cong
Ext^1(I_{\xi}, I_{\xi} \otimes K_S)_0.
\end{eqnarray*}

For simplicity, we assume that $\Supp(\xi) = \{s\} \subset S$.
Note that the vector spaces $Ext^1(I_{\xi}, I_{\xi})_0, 
Ext^1(I_{f^*\xi}, I_{f^*\xi})_0, Ext^2(I_{f^*\xi}, I_{f^*\xi})_0$
all have dimension $2d$.

Applying the local-to-global spectral sequence to
$Ext^1(I_{\xi}, I_{\xi})$, we obtain
\begin{eqnarray*} 
0 \to H^1(S, \mathcal O_S) \to Ext^1(I_{\xi}, I_{\xi})
\to H^0 \big (S, \mathcal Ext^1(I_{\xi}, I_{\xi}) \big ) 
\to H^2(S, \mathcal O_S).
\end{eqnarray*}
It follows that we have an exact sequence
\begin{eqnarray}   \label{obstr.3}
0 \to Ext^1(I_{\xi}, I_{\xi})_0
\to H^0 \big (S, \mathcal Ext^1(I_{\xi}, I_{\xi}) \big ) 
\to H^2(S, \mathcal O_S).
\end{eqnarray}
Since the second term can be computed locally,
by taking $S = \mathbb P^2$, we see that 
\begin{eqnarray*} 
h^0 \big (S, \mathcal Ext^1(I_{\xi}, I_{\xi}) \big ) 
= 2d
\end{eqnarray*}
for an arbitrary surface $S$. So we conclude from 
(\ref{obstr.3}) that
\begin{eqnarray} \label{obstr.4}
Ext^1(I_{\xi}, I_{\xi})_0  \cong
H^0 \big (S, \mathcal Ext^1(I_{\xi}, I_{\xi}) \big )
\end{eqnarray}
since $\dim Ext^1(I_{\xi}, I_{\xi})_0 = 2d$. Similarly, 
we have canonical isomorphisms:
\begin{eqnarray} \label{obstr.5}
       Ext^1(I_{f^*\xi}, I_{f^*\xi})_0  
&\cong&H^0 \big (X, \mathcal Ext^1(I_{f^*\xi}, 
         I_{f^*\xi}) \big )   \nonumber   \\
&\cong&H^0 \big (S, 
         f_*\mathcal Ext^1(I_{f^*\xi}, I_{f^*\xi}) \big ).
\end{eqnarray}

As in (\ref{obstr.3}), we have an injection
\begin{eqnarray*}   
0 \to Ext^1(I_{\xi}, I_{\xi}\otimes K_S)_0
\to H^0 \big (S, \mathcal Ext^1(I_{\xi}, I_{\xi}\otimes K_S) \big ).
\end{eqnarray*}
Note that $H^0 \big (S, \mathcal Ext^1(I_{\xi}, I_{\xi} 
\otimes K_S) \big ) \cong H^0 \big (S, \mathcal Ext^1(I_{\xi}, 
I_{\xi}) \big ) \otimes_\C K_S|_s$ since 
$\mathcal Ext^1(I_{\xi}, I_{\xi})$ is supported
at $\Supp(\xi) = \{s\}$, where $K_S|_s$ is 
the fiber of $K_S$ at $s \in S$. So we get
\begin{eqnarray*}   
0 \to Ext^1(I_{\xi}, I_{\xi}\otimes K_S)_0 \to 
H^0 \big (S, \mathcal Ext^1(I_{\xi}, I_{\xi}) \big )
\otimes_\C K_S|_s.
\end{eqnarray*}
By (\ref{obstr.4}) and the Serre duality, 
$Ext^1(I_{\xi}, I_{\xi}\otimes K_S)_0$ and 
$H^0 \big (S, \mathcal Ext^1(I_{\xi}, I_{\xi}) \big )$
have the same dimension. Hence, we get an isomorphism
\begin{eqnarray}   \label{obstr.6}
Ext^1(I_{\xi}, I_{\xi}\otimes K_S)_0 \cong
H^0 \big (S, \mathcal Ext^1(I_{\xi}, I_{\xi}) \big )
\otimes_\C K_S|_s.
\end{eqnarray}

Again as in (\ref{obstr.3}), we have another injection:
\begin{eqnarray*}  
0 \to Ext^1(I_{f^*\xi}, I_{f^*\xi} \otimes f^*K_S)_0
\to H^0 \big (X, \mathcal Ext^1(I_{f^*\xi}, I_{f^*\xi} 
\otimes f^*K_S) \big ).
\end{eqnarray*}
By the Serre duality, 
$Ext^1(I_{f^*\xi}, I_{f^*\xi} \otimes f^*K_S)_0
\cong Ext^2(I_{f^*\xi}, I_{f^*\xi})_0^{\rm v}$.
Also,
\begin{eqnarray*}  
       H^0 \big (X, \mathcal Ext^1(I_{f^*\xi}, I_{f^*\xi} 
          \otimes f^*K_S) \big )
&\cong&H^0 \big (S, f_*\mathcal Ext^1(I_{f^*\xi}, I_{f^*\xi} 
          \otimes f^*K_S) \big )    \\
&\cong&H^0 \big (S, f_*\mathcal Ext^1(I_{f^*\xi}, I_{f^*\xi}) 
          \otimes K_S\big )    \\
&\cong&H^0 \big (S, f_*\mathcal Ext^1(I_{f^*\xi}, I_{f^*\xi}) 
          \big ) \otimes_\C K_S|_s
\end{eqnarray*}
since $f_*\mathcal Ext^1(I_{f^*\xi}, I_{f^*\xi})$ is 
supported on $\Supp(\xi) = \{s\}$. Therefore, we obtain
\begin{eqnarray}   \label{obstr.7}
0 \to Ext^2(I_{f^*\xi}, I_{f^*\xi})_0^{\rm v}
\to H^0 \big (S, f_*\mathcal Ext^1(I_{f^*\xi}, I_{f^*\xi}) 
\big ) \otimes_\C K_S|_s.
\end{eqnarray}
Since $Ext^2(I_{f^*\xi}, I_{f^*\xi})_0$ and 
$Ext^1(I_{f^*\xi}, I_{f^*\xi})_0$ have the same dimension,
we obtain
\begin{eqnarray*}   
       Ext^2(I_{f^*\xi}, I_{f^*\xi})_0^{\rm v} 
&\cong&H^0 \big (S, f_*\mathcal Ext^1(I_{f^*\xi}, 
         I_{f^*\xi}) \big ) \otimes_\C K_S|_s    \\
&\cong&Ext^1(I_{f^*\xi}, I_{f^*\xi})_0 \otimes_\C K_S|_s
\end{eqnarray*}
from (\ref{obstr.5}) and (\ref{obstr.7}). Combining this with
(\ref{obstr.2}) and (\ref{obstr.4}), we get
\begin{eqnarray}   \label{obstr.8}
Ext^2(I_{f^*\xi}, I_{f^*\xi})_0^{\rm v} \cong
H^0 \big (S, \mathcal Ext^1(I_{\xi}, I_{\xi}) \big )
\otimes_\C K_S|_s.
\end{eqnarray}
In view of (\ref{obstr.6}), the Serre duality and 
(\ref{obstr.2}), we conclude that
\begin{eqnarray*}
      Ext^2(I_{f^*\xi}, I_{f^*\xi})_0
\cong Ext^1(I_{\xi}, I_{\xi}\otimes K_S)_0^{\rm v}  
\cong Ext^1(I_{\xi}, I_{\xi})_0
\cong Ext^1(f^*I_{\xi}, f^*I_{\xi})_0. 
\end{eqnarray*}
This completes the proof of the isomorphism (\ref{obstr.1}).
\end{proof}

Next, we consider the case when either $n \ne 0$ or 
$\beta \ne d \beta_0$ with $d \ge 0$.
We further assume that the moduli space $\mathfrak I_n(X, \beta)$ 
is nonempty. For simplicity, put
\begin{eqnarray*}  
\mathfrak I = \mathfrak I_n(X, \beta).
\end{eqnarray*}
Let $\mathcal I$ be the universal ideal sheaf over 
$\mathfrak I \times X$. Denote the trace-free part of 
the element $R \mathcal Hom(\mathcal I,\mathcal I)$ 
in the derived category $\mathfrak D_{\rm coh}
(\mathfrak I \times X)$ by 
\begin{eqnarray*}
R \mathcal Hom(\mathcal I,\mathcal I)_0.
\end{eqnarray*}
Let $\pi: \mathfrak I \times X \to \mathfrak I$ be 
the projection. By \cite{Tho}, the virtual fundamental
class $[\mathfrak I]^{\rm vir}$ is defined via 
the following element in the derived category 
$\mathfrak D_{\rm coh}(\mathfrak I)$:
\begin{eqnarray}  \label{E}
  \mathcal E 
= R\pi_* \big ( 
  R\mathcal Hom(\mathcal I,\mathcal I)_0 \big ).
\end{eqnarray}

Let $p \in E$, and consider the sheaf $\big ( \id_\mathfrak I 
\times \phi_p \times \id_S \big )^*\mathcal I$ over 
\begin{eqnarray*}
\mathfrak I \times X = \mathfrak I \times E \times S.
\end{eqnarray*}
We see from (\ref{fix_beta}) that $\big ( \id_\mathfrak I 
\times \phi_p \times \id_S \big )^*\mathcal I$
is a flat family of ideal sheaves whose corresponding
$1$-dimensional closed subschemes satisfy (\ref{In}). 
By the universal property of the moduli space 
$\mathfrak I$, there is an automorphism 
\begin{eqnarray}  \label{Phip}
\Phi_p: \mathfrak I \to \mathfrak I
\end{eqnarray}
such that $(\Phi_p \times \id_X)^*\mathcal I 
= \big ( \id_\mathfrak I \times \phi_p \times \id_S 
\big )^*\mathcal I \cong \mathcal I$. In particular, 
$E$ acts on $\mathfrak I$.

\begin{lemma} \label{lma_dt_0}
Let $n \ne 0$ or $\beta \ne d \beta_0$ with $d \ge 0$. Then,
\begin{eqnarray}   \label{lma_dt_0.1}
\langle \wtau_{k_1}(\gamma_1) \cdots 
   \wtau_{k_r}(\gamma_r) \rangle_{n, \beta}
= 0
\end{eqnarray}
whenever $\gamma_1, \ldots, \gamma_r \in f^*H^*(S; \Q)
\subset H^*(X; \Q)$.
\end{lemma}
\begin{proof}
The proof is similar to that of Lemma~\ref{lma_gw=0}.
Assume that the moduli space 
\begin{eqnarray*}
\mathfrak I = \mathfrak I_n(X, \beta)
\end{eqnarray*}
is nonempty. 
If $\beta \ne d \beta_0$ with $d \ge 0$, then the algebraic
group $E$ acts on $\mathfrak I$ with finite stabilizers. 
If $\beta = d \beta_0$ with 
$d \ge 0$ and if $I_Z \in \mathfrak I$, then $Z$ consists
of a curve $f^*(\xi)$ for some $\xi \in S^{[d]}$ 
and of some (possibly embedded) points of length $n \ne 0$.
So again $E$ acts on the moduli space $\mathfrak I$ 
with finite stabilizers.

As in the proof of Lemma~\ref{lma_gw=0}, there exists
some number $N$ such that if $G$ is a cyclic subgroup of 
$E$ of prime order $p > N$, then $G$ acts freely on 
$\mathfrak I$. Fix such cyclic subgroups $G$ of $E$.
Since the complex $R\pi_* \big ( R\mathcal Hom
(\mathcal I,\mathcal I)_0 \big )$ from 
(\ref{E}) is equivariant for the action of 
any algebraic automorphism group of $X$, 
the cycle $[\mathfrak I]^{\rm vir}$
defines an element of the equivariant Borel-Moore 
homology group $H_*^G(\mathfrak I)$. For $1 \le i \le r$,
choose a positive integer $m_i$ such that the multiple 
$m_i \gamma_i$ defines an element of $H^*_G(X)$.
It follows from (\ref{def_opr2}) and (\ref{dt_desc}) 
that the cycle 
\begin{eqnarray*}
m_1 \cdots m_r \; \langle \wtau_{k_1}(\gamma_1) 
\cdots \wtau_{k_r}(\gamma_r) \rangle_{n, \beta}
\end{eqnarray*}
defines an element in the degree-$0$ Borel-Moore
homology $H^G_0(\mathfrak I)$. 
Again as in the proof of Lemma~\ref{lma_gw=0},
we conclude that $\langle \wtau_{k_1}(\gamma_1) 
\cdots \wtau_{k_r}(\gamma_r) \rangle_{n, \beta} = 0$.
\end{proof}

\begin{proposition} \label{prop_dt2}
Let $\beta \in H_2(X; \Z) \backslash \{0\}$. Assume 
$\displaystyle{\int_{\beta} K_X} 
= \displaystyle{\int_{\beta} f^*K_S} =0$. Then,
\begin{eqnarray*}
&&  {\bf Z}_{\text{\rm DT}}'\left (X; q| 
   \prod_{i=1}^r \wtau_{0}(\gamma_i) \right )_\beta \\
&=&\left \{
   \begin{array}{ll}
     \prod_{i=1}^r \int_{\beta} \gamma_i 
       \cdot \chi \big ( S^{[d]} \big )
       &\text{if} \; |\gamma_i| = 2 \; 
        \text{for every} \; i \text{ and } \beta = d\beta_0
     \text{ for some } d \ge 1;\\
     0 &\text{otherwise}.
   \end{array}
  \right.
\end{eqnarray*}
\end{proposition}
\begin{proof}
First of all, since $\chi(X) = 0$, we see from (\ref{rpf_dt}) 
and (\ref{p_func}) that
\begin{eqnarray}     \label{prop_dt2.1}
  {\bf Z}_{\text{\rm DT}}'\left (X; q| 
   \prod_{i=1}^r \wtau_{k_i}(\gamma_i) \right )_\beta
= \sum_{n \in \Z} \left \langle \wtau_{k_1}(\gamma_1) \cdots 
   \wtau_{k_r}(\gamma_r) \right \rangle_{n, \beta} \; q^n.
\end{eqnarray}

Next, in view of (\ref{dim_In}) and the condition on degrees, 
we have 
\begin{eqnarray*}
\sum_{i=1}^r |\gamma_i| = 2r, \qquad |\gamma_r| \le 2.
\end{eqnarray*}
Therefore, we conclude from (\ref{dt_desc}) and 
Lemma~\ref{dt_fund_div}~(i) that
\begin{eqnarray*}  
\langle \wtau_{0}(\gamma_1) \cdots 
\wtau_{0}(\gamma_r) \rangle_{n, \beta} 
= \left \{
   \begin{array}{ll}
     \prod_{i=1}^r \int_{\beta} \gamma_i
       \cdot \langle  \rangle_{n, \beta}
       &\text{\rm if} \; |\gamma_i| = 2 \; 
        \text{\rm for every} \; i;\\
     0 &\text{\rm otherwise}.
   \end{array}
  \right.
\end{eqnarray*}
By Lemma~\ref{obstr}~(ii) and Lemma~\ref{lma_dt_0},
we obtain
\begin{eqnarray*}  
&&\langle \wtau_{0}(\gamma_1) \cdots 
  \wtau_{0}(\gamma_r) \rangle_{n, \beta}    \\
&=&\left \{
   \begin{array}{ll}
     \prod_{i=1}^r \int_{\beta} \gamma_i
       \cdot \chi \big ( S^{[d]} \big )
       &\text{\rm if} \; |\gamma_i| = 2 \; 
        \text{\rm for every} \; i, n=0, \beta = d\beta_0;\\
     0 &\text{\rm otherwise}.
   \end{array}
  \right.
\end{eqnarray*}
Now the proposition follows immediately from (\ref{prop_dt2.1}).
\end{proof}

\begin{proposition} \label{prop_dt1}
Let $X$ be from Notation~\ref{nota} and 
$\beta \in H_2(X; \Z) \backslash \{0\}$.
Assume that $\gamma_1, \ldots, \gamma_r \in f^*H^*(S; \Q)
\subset H^*(X; \Q)$. Then,
\begin{eqnarray*}
  {\bf Z}_{\text{\rm DT}}'\left (X; q| 
   \prod_{i=1}^r \wtau_{k_i}(\gamma_i) \right )_\beta
= \left \{
   \begin{array}{ll}
     \chi \big ( S^{[d]} \big ) &if \; r = 0 \text{ and }
       \beta = d\beta_0 \text{ with } d \ge 1,\\
     0 &otherwise.
   \end{array}
  \right.
\end{eqnarray*}
\end{proposition}
\begin{proof}
If $\beta \ne d\beta_0$, then the proposition follows from
(\ref{prop_dt2.1}) and Lemma~\ref{lma_dt_0}. In the rest of
the proof, we let $\beta = d\beta_0$ with $d \ge 1$.
By (\ref{prop_dt2.1}) and Lemma~\ref{lma_dt_0} again,
\begin{eqnarray*}   
  {\bf Z}_{\text{\rm DT}}'\left (X; q| 
   \prod_{i=1}^r \wtau_{k_i}(\gamma_i) \right )_{d\beta_0}
= \left \langle \wtau_{k_1}(\gamma_1) \cdots 
   \wtau_{k_r}(\gamma_r) \right \rangle_{0, d\beta_0}.
\end{eqnarray*}
Thus we see from Lemma~\ref{obstr}~(ii) that the proposition
holds if $r = 0$.

To prove our proposition, it remains to verify that 
if $r \ge 1$, then
\begin{eqnarray}   \label{prop_dt1.1}
\left \langle \wtau_{k_1}(\gamma_1) \cdots 
\wtau_{k_r}(\gamma_r) \right \rangle_{0, d\beta_0} = 0.
\end{eqnarray}
Since the expected dimension of $\mathfrak I_0(X, d\beta_0)$ 
is zero, (\ref{prop_dt1.1}) holds unless
\begin{eqnarray}  \label{prop_dt1.2}
\sum_{i=1}^r (2k_i - 2 + |\gamma_i|) = 0, \qquad
(2k_r - 2 + |\gamma_r|) \le 0.
\end{eqnarray}
W.l.o.g., we may assume that $k_{\w r +1} = k_{\w r +2} =
\ldots = k_{r} = 0$
and
\begin{eqnarray}   \label{prop_dt1.3}
k_1, \ldots, k_{\w r - 1}, k_{\w r} \ge 1
\end{eqnarray}
for some $\w r$ with $0 \le \w r \le r$. Then we see from 
(\ref{dt_desc}), Lemma~\ref{dt_fund_div}~(i) and 
(\ref{prop_dt1.2}) that (\ref{prop_dt1.1}) holds unless
$\w r = r$, $k_1 = \ldots = k_r= 1$, and
$|\gamma_1| = \ldots = |\gamma_r| = 0$.
When 
\begin{eqnarray*}
k_1 = \ldots = k_r= 1, 
\quad |\gamma_1| = \ldots = |\gamma_r| = 0, 
\end{eqnarray*}
(\ref{prop_dt1.1}) follows from Lemma~\ref{ch3}~(ii) since 
the moduli space $\mathfrak I_0(X, d\beta_0)$ is smooth.
\end{proof}

\end{document}